\documentclass{amsart}

\usepackage{amsfonts}
\usepackage{latexsym}

\newtheorem{theorem}{Theorem}[section]
\newtheorem{lemma}[theorem]{Lemma}

\newtheorem{proposition}[theorem]{Proposition}

\theoremstyle{definition}

\theoremstyle{remark}
\newtheorem{remark}[theorem]{Remark}

\numberwithin{equation}{section}

\DeclareMathOperator{\dist}{dist}

\DeclareMathOperator{\Ran}{Ran}
\DeclareMathOperator{\rank}{rank}
\DeclareMathOperator{\Rank}{Rank}
\DeclareMathOperator{\Tr}{Tr}

\begin{document}

\title[Scattering Poles for Hyperbolic Manifolds]{Scattering Poles for Asymptotically Hyperbolic Manifolds}
\author[D. Borthwick]{David Borthwick}
\address{Department of Mathematics and Computer Science,
Emory University, Atlanta, Georgia 30322, U. S. A.}
\email{davidb@mathcs.emory.edu}
\thanks{Supported in part by NSF grant DMS-9796195 and by an NSF
Postdoctoral Fellowship.}
\author{Peter Perry}
\thanks{Supported in part by NSF grant DMS-9707051}
\address{Department of Mathematics,
	University of Kentucky,
	Lexington, Kentucky 40506--0027,
	U. S. A}
\email{perry@ms.uky.edu}
\date{May 25, 1999}

\begin{abstract}
For a class of manifolds $X$ that includes quotients of real hyperbolic $
(n+1)$-dimensional space by a convex co-compact discrete group, we show that
the resonances of the meromorphically continued resolvent kernel for the
Laplacian on $X\,$ coincide, with multiplicities, with the poles of the
meromorphically continued scattering operator for $X$. In order to carry out
the proof, we use Shmuel Agmon's perturbation theory of resonances to show
that both resolvent resonances and scattering poles are simple for generic
potential perturbations.
\end{abstract}

\maketitle

\section{Introduction}

The purpose of this paper is to show the equivalence of two possible
notions of `scattering resonances' for the Laplacian on asymptotically
hyperbolic manifolds, i.e., complete Riemannian manifolds of infinite
volume with `constant curvature at infinity'. On the one hand, it
is very natural to define scattering resonances with respect to
the meromorphically continued resolvent of the Laplace operator.
This point of view has been very fruitful and has led to a large
body of results on the distribution of scattering resonances in
the complex plane;  see \cite{Zworski:1994} for a survey and
\cite{GZ:1995a,GZ:1995b,GZ:1997,Zworski:1997} for results on the
class of manifolds studied here. On the other hand, it is also
reasonable, by analogy with scattering theory for the Schr\"{o}dinger
or wave equations in Euclidean space (see e.g. \cite{LP:1967,RS:1979}),
to define scattering resonances as poles of the scattering operator
for the Laplacian.  For Schr\"odinger scattering on Euclidean space,
the equivalence of scattering resonances and resolvent resonances
is well-known \cite{Hagedorn:1979,Jensen:1977,Jensen:1980}.  The
analogous result for hyperbolic manifolds is of interest since the
poles of the scattering operator have a geometric and dynamical interpretation:
they are among the poles of Selberg's zeta function for geodesic
flow on the manifold \cite{Patterson:1991,PP:1998,Perry:1990}. Thus the
scattering resonances serve, in a sense, as discrete data similar
in character to the eigenvalues of a compact surface.

For non-compact Riemann surfaces and certain metric perturbations,
Guillop\'{e} and Zworski showed that the set of resolvent resonances
and the set of scattering poles coincide with multiplicities
(\cite{GZ:1997}, Proposition 2.11); here we will show that, with
some restrictions, the same
equivalence holds for asymptotically hyperbolic manifolds in higher
dimension.

To describe our results, we first recall that an asymptotically hyperbolic
manifold is a compact manifold $\overline{X}$ with boundary endowed with a
Riemannian metric of a special form. A {\em defining function} for the
boundary of a compact manifold $\overline{X}$ is a nonnegative 
$C^{\infty }$ function on $\overline{X}$ with 
$x^{-1}(0)=\partial \overline{X}$ and 
$\left.  dx \right| _{\partial \overline{X}}$ nowhere vanishing. 
The metric $g$ then takes the form $x^{-2}h$ where 
$x$ is a defining function and $h$ is a nondegenerate smooth
metric on $\overline{X}$ such that $\left| dx \right| _{h}\rightarrow 1$
as $x \rightarrow 0$. Note that this metric puts $\partial \overline{X}$
`at infinity' and makes $X$, the interior of $\overline{X}$, a complete
Riemannian manifold of infinite volume. The condition on
$\left|dx\right|_h$ insures
that the sectional curvatures approach $-1$ at metric infinity. If $\Delta
_{g}$ denotes the positive Laplace-Beltrami operator on $(X,g)$ and $X$ has
dimension $n+1$, it is known that the spectrum of $\Delta _{g}$ consists of
at most finitely many $L^{2}$ eigenvalues of finite multiplicity in the
interval $[0,n^{2}/4)$ (see \cite{LP:1982} for quotients of hyperbolic 
space and \cite{Mazzeo:1991b} for asymptotically hyperbolic manifolds) and
absolutely continuous spectrum in $[n^{2}/4,\infty )$ (see 
\cite{FHP:1991,LP:1984}
for hyperbolic quotients and \cite{MM:1987} for resolvent 
estimates which imply absolute continuity of the spectrum for 
asymptotically hyperbolic manifolds). Thus the resolvent 
${\mathcal R}_{g}(z)=(\Delta _{g}-z)^{-1}$ is a meromorphic 
operator-valued function on the cut plane 
${\mathbb C}\backslash \lbrack n^{2}/4,\infty )$ with poles at the $L^{2}$
eigenvalues having finite-rank residues. It is convenient to introduce a
uniformizing parameter $\zeta $ and set $R_{g}(\zeta )=(\Delta _{g}-\zeta
(n-\zeta ))^{-1}$, which is then a meromorphic operator-valued function on
the half-plane $\Re (\zeta )>n/2$. The operator $R_{g}(\zeta )$ has
first-order poles at points $\zeta _{0}$ whenever 
$\zeta _{0}(n-\zeta _{0})$
is an $L^{2}(X)$-eigenvalue of $\Delta _{g}$. We denote by $Z_{p}$ the
(finite and possibly empty) set of all such $\zeta _{0}$. The
\emph{multiplicity} of $\zeta _{0}\in Z_{p}$ is the dimension, 
$m_{\zeta _{0}}$, of the eigenspace of $\Delta _{g}$ with eigenvalue 
$\zeta _{0}(n-\zeta _{0})$. Equivalently,
\begin{equation}
\label{eq.r.ev.mult}
m_{\zeta _{0}}={{\rank}}
	\left( \int_{\gamma _{\zeta _{0}}}\, \,R_{g}(\zeta )\,d\zeta \right), 
\end{equation}
where $\gamma_{\zeta_0}$ is a simple closed contour surrounding $\zeta_0$
and no other pole of $R_g(\zeta)$.

First, we define the resolvent resonance set of $\Delta _{g}$. 
Let
$\dot{C}^{\infty }(X)$ denote
the smooth functions on $\overline{X}$ vanishing to all orders at $\partial
\overline{X} $. Viewed as a map from 
$\dot{C}^{\infty }(X)$ to $C^{\infty }(X)$, the
resolvent operator $R_{g}(\zeta )$ admits a meromorphic continuation to 
${\mathbb C}{\bf -}\frac{1}{2}(n-{\mathbb N})$, as was shown in 
\cite{MM:1987}.
Resolvent resonances are poles of this meromorphic continuation in the
half-plane $\Re (\zeta )<n/2$, excluding the region
$\frac{1}{2}(n-{\mathbb N})$. If $X$ has constant curvature in a
neighborhood of infinity, the analysis of \cite{GZ:1995b} shows that the 
resolvent has a meromorphic continuation to $\Re(\zeta)<n/2$ with no 
further restriction. Lower bounds on resolvent resonances proven
in \cite{SJ:1993} show that this set is always nontrivial for constant
curvature spaces, and explicit examples (see, for instance, section 3 of \cite
{GZ:1995a}) show that resolvent resonances can form an infinite lattice in
the half-plane $\Re (\zeta )<n/2$. If $\zeta _{0}$ is a resolvent resonance
and $\Re (\zeta _{0})<n/2$, the {\em multiplicity} of $\zeta _{0}$
is the number 
\begin{equation}
m_{\zeta _{0}}={{\rank}}\left( \int_{\gamma _{\zeta _{0}}}\,
\,R_{g}(\zeta )\,d\zeta \right)  \label{eq.r.mult}
\end{equation}
(cf. (\ref{eq.r.ev.mult})). Here $\gamma _{\zeta _{0}}$ is a simple closed
curve that encloses $\zeta _{0}$ and no other pole of $R_{g}(\zeta )$. 
The point $\zeta_0$ is a {\em semi-simple resonance} if 
$R_{g}(\zeta)$ has a simple pole at $\zeta _{0}$. The point
$\zeta_0$ is a {\em simple resonance} if, in addition,  
the residue of the pole has rank one.
The {\em resolvent resonance set} is the set ${\mathcal R}$ of 
all $\zeta _{0}$ such that $\Re (\zeta _{0})<n/2$, $\zeta_0 \not\in
\frac{1}{2}(n-{\mathbb N})$,  and $R_{g}(\zeta )$ has 
a pole of multiplicity $m_{\zeta_{0}}\neq 0$ at $\zeta _{0}$.

Next, we define the scattering operator and the scattering resonance 
set of $\Delta_{g}$. For $\zeta \in {\mathbb C}$ with $\Re (\zeta )=n/2$ 
and $\zeta \neq n/2$, and each 
$f_{-}\in C^{\infty }(\partial \overline{X})$, there is
a unique smooth solution of the eigenvalue equation $(\Delta _{g}-\zeta
(n-\zeta ))u=0$ having the asymptotic form 
$$
u=x^{n-\zeta }f_{+}+x ^{\zeta }f_{-}+O(x ^{n/2+1}) 
$$
where $f_{+}\in C^{\infty }(\partial \overline{X})$; for a proof see \cite
{Borthwick:1997} or \cite{JS:1998}. It follows that $f_{+}$ is uniquely
determined and that there is a linear map $S(\zeta ):C^{\infty }(\partial 
\overline{X})\rightarrow C^{\infty }(\partial \overline{X})$ with $S(\zeta
)f_{-}=f_{+}$; moreover it is clear that $S(\zeta )S(n-\zeta )=I$. It can 
be shown that $S(\zeta )$ extends to a meromorphic family of operators on 
${\mathbb C}$
(see \cite{Borthwick:1997} or \cite{JS:1998}); these operators may have
infinite-rank poles at $\zeta \in n/2+{\mathbb N}$, and infinite-rank 
zeros at $\zeta \in n/2-{\mathbb N}$. A scattering pole is a pole of the 
meromorphic continuation of $S(\zeta )$ in the half-plane 
$\Re(\zeta )<n/2$, excluding the set $\frac{1}{2}(n-{\mathbb N})$ unless
$X$ has constant curvature near infinity. It can be
shown that the continued operator admits the factorization 
\begin{equation}
S(\zeta )=P(\zeta )(I+K(\zeta ))Q(\zeta )  \label{eq.s.factor}
\end{equation}
holds, where $P(\zeta )$ and $Q(\zeta )$ are holomorphically invertible
families of elliptic operators for 
$\zeta \in {\mathbb C}-\frac{1}{2}(n-{\mathbb N})$ and $K(\zeta )$ is a 
meromorphic family of compact operators on 
$C^{\infty }(\partial \overline{X})$. For 
$\zeta _{0}\notin \frac{1}{2}(n-{\mathbb N})$, the {\em multiplicity} of 
a scattering pole $\zeta _{0}$ is the integer 
\begin{equation}
\nu_{\zeta _{0}}=\frac{1}{2\pi i}\,{{\Tr}}\left( \int_{\gamma
_{\zeta _{0},\varepsilon }}\,S(\zeta )^{-1}S^{\prime }(\zeta )\,d\zeta
\right)  \label{eq.s.mult}
\end{equation}
(compare \cite{GZ:1997,PP:1998} where similar definitions are made).
The factorization (\ref{eq.s.factor}) and results of \cite{GS:1972}
show that $\nu _{\zeta _{0}}$ is an integer equal to the multiplicity
of zeros minus the multiplicity of poles of $(I+K(\zeta ))$ at
$\zeta =\zeta _{0}$ (we give a precise formulation in section
\ref{sec.scattering}).  We will say that $\zeta _{0}$ is semi-simple
if the pole of $S(\zeta )$ is of first order, and simple if the
residue is rank-one. The {\em scattering resonance set} is the set
${\mathcal S}$ of all pairs $(\zeta _{0},-\nu_{\zeta_0})$ with
$\Re(\zeta_0)< n/2$, $\zeta \not\in \frac{1}{2}(n-{\mathbb N})$, 
and $\nu _{\zeta _{0}}\neq 0$.

We would like to show a correspondence, with multiplicities, between the
sets ${\mathcal R}$ and ${\mathcal S}$. A direct method (see for example 
\cite{GZ:1997}, where the case $n=1$ is treated) would compare the 
Laurent expansion of the meromorphically continued resolvent at 
$\zeta _{0}\in {\mathcal R}$ to the Laurent expansion of the scattering 
operator at $\zeta _{0}$, using the fact that the Schwarz kernel of the 
scattering operator can be recovered from that of the resolvent kernel. 
This direct method works easily when the resolvent resonance is simple but is 
somewhat complicated for non-simple resonances. For this reason, we will 
perturb the operator $\Delta _{g}$ with a potential 
$V\in \dot{C}^{\infty }(X)$ which, as we will show, can
be chosen to make all resonances of the meromorphically continued resolvent
$$
R_{V}(\zeta )=(\Delta _{g}+V-\zeta (n-\zeta ))^{-1}
$$
simple. The perturbation will split each resonance of multiplicity $m$  into 
$m$ resonances of multiplicity one localized near the unperturbed resonance,
and similarly each eigenvalue of multiplicity $m$ into $m$ eigenvalues of
multiplicity one. This result, which is of some independent interest, will
allow us to count multiplicities properly but avoid technicalities
associated with non-simple resonances.

Our analysis of generic potential perturbations is inspired by Klopp and
Zworski's analysis of resonances in potential scattering \cite{KZ:1995}. To
carry out the analysis, we will use Shmuel Agmon's perturbation theory of
resonances \cite{Agmon:1997} in which the resonances are realized as
eigenvalues of a non-self-adjoint operator on a cleverly constructed Banach
space; this replaces the complex scaling used in \cite{KZ:1995}. 
Standard Kato-Rellich perturbation theory \cite{Kato:1976} can then be
used to study how the resonances move under perturbation.

Our first result is:

\begin{theorem}
\label{thm.essential}
Let $(X,g)\,$ be an asymptotically hyperbolic manifold,
and let ${\mathcal R}$ and ${\mathcal S}$ be respectively the resolvent 
resonance set and scattering resonance set for the Laplacian $\Delta _{g}$.
Then ${\mathcal R}={\mathcal S}$ on 
${\mathbb C}-\frac{1}{2}(n-{\mathbb N})$ except for at most finitely
many points. More precisely,
for any $\zeta _{0}$ with $\Re (\zeta _{0})<n/2$, 
$\zeta _{0}\notin \frac{1}{2}(n-{\mathbb N})$, 
the relationship 
$$
\nu_{\zeta _{0}}=m_{n-\zeta _{0}}-m_{\zeta _{0}} 
$$
holds.
\end{theorem}

Note that $m_{n-\zeta _{0}}$ is nonzero only for the finitely many
$\zeta _{0}$ with $n-\zeta _{0}\in Z_{p}$; for all other $\zeta
_{0}\in \frac{1}{2}(n-{\mathbb N})$, the scattering resonances and
resolvent resonances coincide with multiplicities.

We can make a stronger statement if $(X,g)$ has even dimension and constant
curvature in a neighborhood of infinity, i.e., if there is a compact subset 
$K$ of $X$ so that $g$ has constant negative curvature $-1$ on 
$X\backslash K$. This class was studied in \cite{GZ:1995b} and
includes convex co-compact hyperbolic 
manifolds. Recall that a geometrically finite group $\Gamma $ of 
isometries of real hyperbolic $(n+1)$--dimensional space 
${\mathbb H}^{n+1}$ is called convex co-compact if the orbit space 
$\Gamma \backslash {\mathbb H}^{n+1}$ has infinite volume and 
$\Gamma$ contains no parabolic elements. If $\Gamma $ is torsion-free 
(which we can insure by passing to a subgroup of finite index), the orbit 
space $X=\Gamma \backslash {\mathbb H}^{n+1}$ is a complete Riemannian 
manifold when given the induced hyperbolic metric $g$. 

To formulate a result, we introduce a renormalized scattering operator 
\[
S_r(\zeta )=\frac{\Gamma (s-n/2)}{\Gamma (n/2-s)}S(\zeta ); 
\]
if $\dim(X)$ is even (so $n$ is odd), the renormalization has the effect 
of dividing out the 
infinite-rank zeros of $S(\zeta )$ at $\zeta \in n/2-{\mathbb N}$. 
Note, however, that
it may also cancel poles of the scattering operator due to eigenvalues of
$\Delta_g$ if the set $Z_p \cap (\frac{1}{2}n+{\mathbb N})$ is nonempty.
It can be shown that the factorization 
\[
S_r(\zeta )=P(\zeta )(I+K(\zeta ))P(\zeta ) 
\]
holds for a family of operators $P(\zeta )$ holomorphically invertible in
all of ${\mathbb C}$, and a meromorphic family of compact operators $K(\zeta )$
with at most finite rank zeros and poles. We then define the multiplicity of
a scattering pole using (\ref{eq.s.mult}) but with $S(\zeta )$ replaced by $
S_r(\zeta )$ in the definition, and enlarge the set ${\mathcal S}$ to include 
any $\zeta _{0}$ with $\Re (\zeta _{0})<n/2$ where $S_r(\zeta )$ has a
pole. Similarly, we enlarge the set ${\mathcal R}$ to include all resolvent poles
$\zeta$ with $\Re(\zeta) < n/2$. With these definitions, we have:

\begin{theorem}
\label{thm.hyperbolic}
Let $(X,g)$ have constant curvature near infinity and
suppose that $\dim (X)$ is even. Let ${\mathcal R}$ and ${\mathcal S}$ denote the
resolvent resonance and scattering resonance sets for $\Delta _{g}$. Then
the relation 
$$
\nu _{\zeta _{0}}=m_{n-\zeta _{0}}-m_{\zeta _{0}} 
$$
holds for all $\zeta _{0}$ with 
$n-\zeta _{0}\notin Z_{p}\cap (\frac{1}{2}n+ {\mathbb N})$.
\end{theorem}

Thus, for all but finitely many $\zeta _{0}$, the resolvent resonance set
and the scattering resonance set coincide with multiplicities. The set $
Z_{p}\cap \left( \frac{1}{2}n+{\mathbb N}\right) $ consists at most of finitely
many elements, and is empty if $n=1$. 

This paper is organized as follows. In section \ref{sec.resolvent} we review
the Mazzeo-Melrose construction of the resolvent and study its behavior near
resolvent resonances. In section \ref{sec.scattering} we recall how the
scattering operator can be recovered from the resolvent and discuss its
behavior near scattering poles. In section \ref{sec.perturb} we study the
perturbation behavior of resonances when the operator $\Delta _{g}$ is
perturbed by a potential $V\in \dot{C}^{\infty }(X)$. In section \ref
{sec.generic}, we show that the operator $P_V=\Delta _{g}+V$ has only simple
resonances $\zeta _{0}$ for $\zeta _{0}\in \frac{1}{2}(n-{\mathbb N})$ for
potentials $V$ in a dense open subset of $\dot{C}^{\infty }(X)$. Finally, in
section \ref{sec.equality}, we prove Theorems \ref{thm.essential} and \ref
{thm.hyperbolic}.

In what follows, $x ^{N}L^{2}(X)$ denotes the space of locally
square-integrable functions $v$ on $X$ with $v=x ^{N}u$ for a function $
u\in L^{2}(X)$ and a fixed real number $N$. For a fixed, given $N$, we denote 
by $B_{0}$ the Banach
space $x ^{N}L^{2}(X)$, and by $B_{1}$ the Banach space $x
^{-N}L^{2}(X)$. If ${\mathcal X}$ and ${\mathcal Y}$ are Banach spaces, ${\mathcal L}(
{\mathcal X},{\mathcal Y})$ denotes the Banach space of bounded operators from $
{\mathcal X}$ to ${\mathcal Y}$. We metrize $\dot{C}^{\infty }(X)$ by introducing
the seminorms 
\begin{equation}
\label{eq.seminorm}
d_{\alpha ,n}(V)=\sup_{x\in X}\left| x ^{-n}D^{\alpha }V(x)\right| 
\end{equation}
for nonnegative integers $n$ and multi-indices $\alpha $, and we denote by $
d(\,\cdot \,,\,\cdot \,)$ the associated metric on $\dot{C}^{\infty }(X)$.

\section{The Resolvent of $\Delta _g+V$ and its Meromorphic Continuation}

\label{sec.resolvent}

The resolvent of the operator $P_V=\Delta _{g}+V$ has a distribution kernel,
with respect to the Riemannian density on $X$, which is smooth away from the
diagonal $\Lambda $ of $X\times X$. To describe its singularities, it is
useful to introduce the blow up of $\overline{X}\times \overline{X}$ along 
$ \partial \overline{X}\times \partial \overline{X}$, the stretched product $
\overline{X}\times _{0}\overline{X}$. This amounts to introducing polar
coordinates at the diagonal in the corner of $\overline{X}\times \overline{X}
$ where $\Lambda $ intersects the `top' boundary face $\overline{X}\times
\partial \overline{X}$ and the `bottom' boundary face $\partial \overline{X}
\times \overline{X}$; globally one replaces $\partial \Lambda $ with the
doubly inward-pointing spherical normal bundle of $\partial \Lambda $. If $
(x,y)$ and $(x^{\prime },y^{\prime })$ are local coordinates on $\overline{X}
$ in a neighborhood of the boundary, $\partial \Lambda $ is given by $
x=x^{\prime }=y-y^{\prime }=0$; local coordinates for $\overline{X}\times
_{0}\overline{X}$ near the boundary are then given by $(r,\eta ,\eta
^{\prime },\theta ,y)$ where 
\[
\begin{array}{c}
r=\sqrt{x^{2}+(x^{\prime })^{2}+\left| y-y^{\prime }\right| ^{2}} \\ 
\\ 
(\eta ,\eta ^{\prime },\theta )=(x/r,x^{\prime }/r,(y-y^{\prime })/r).
\end{array}
\]
We denote by $\beta $ the `blow-down map' $\beta :\overline{X}\times _{0}
\overline{X}\rightarrow \overline{X}\times \overline{X}$; in the local
coordinates described above, $\beta (r,\eta ,\eta ^{\prime },\theta
,y)=(r\eta ,y,r\eta ^{\prime },y-r\theta )$.

The following theorem summarizes the Mazzeo-Melrose \cite{MM:1987} construction of the
resolvent. Although Mazzeo and Melrsoe did not treat potential
perturbations, potentials in the class $\dot{C}^{\infty }(X)$ may be
accommodated without difficulty (see \cite{JS:1998}, Theorem 3.1 and its
proof). We denote by $G_{\zeta }$ the integral kernel of the resolvent
operator $R_{V}(\zeta )=(P_V-\zeta (n-\zeta ))^{-1}$ with respect
to Riemannian measure on $X$, initially a meromorphic function of $\zeta $
with $\Re (\zeta )>n/2$.

\begin{theorem}
\label{thm.resolvent.ac}Let $(X,g)$ be an asymptotically hyperbolic
manifold, and let $V\in \dot{C}^{\infty }(X)$. The resolvent kernel $
G_{\zeta }$ has a meromorphic continuation to ${\mathbb C}$ with 
\[
\beta ^{\ast }G_{\zeta }=A_{\zeta }+B_{\zeta }+C_{\zeta } 
\]
where 
\[
A_{\zeta }\in I^{-2}(\overline{X}\times _{0}\overline{X}), 
\]
\[
B_{\zeta }\in (\eta \eta ^{\prime })^{\zeta }C^{\infty }(\overline{X}\times
_{0}\overline{X}), 
\]
and 
\[
C_{\zeta }\in \beta ^{\ast }\left[ (xx^{\prime })^{\zeta }C^{\infty }(
\overline{X}\times \overline{X})\right] . 
\]
Moreover $A_{\zeta }$ is an entire function of $\zeta $, $B_{\zeta }$ is
holomorphic in ${\mathbb C}{\bf -}\frac{1}{2}(n-{\mathbb N})$, and $C_{\zeta }$ is
meromorphic in ${\mathbb C}-\frac{1}{2}(n-{\mathbb N})$.
\end{theorem}

{\em Sketch of proof:} Let $P_{\zeta }=\Delta _{g}+V-\zeta (n-\zeta )$.
Given an operator $A$ on $C^{\infty }(X)$, we will denote by $\kappa (A)$
the lift of the kernel of $A$ (with respect to the Riemannian density on $X$
) to $X\times _{0}X$. The construction in \cite{MM:1987} may be broken into
three pieces. First, we construct an operator $A_{\zeta }$ to cancel the
conormal singularity of $\kappa (P_{\zeta })$ on the lifted diagonal. The
family $A_{\zeta }$ is entire and has the property that $\kappa (I-P_{\zeta
}A_{\zeta })\in C^{\infty }(\overline{X}\times _{0}\overline{X})$. This
remainder does not yet correspond to the integral kernel of a compact
operator on the original space.

To improve the error term, one uses the model resolvent. A second operator $
B_{\zeta }$ is constructed so that $E_{\zeta }=I-P_{\zeta }(A_{\zeta
}+B_{\zeta })$ has 
\[
\kappa (E_{\zeta })\in \eta ^{\zeta }(\eta ^{\prime }r)^{\infty }C^{\infty }(
\overline{X}\times _{0}\overline{X}). 
\]
The operator $B_{\zeta }$ is holomorphic in ${\mathbb C}-\frac{1}{2}(n-{\mathbb N}
) $; the operator $E_{\zeta }$ is a compact operator on the weighted $L^{2}$
space $x ^{N}L^{2}(X)$ for all $\zeta $ with $\Re (\zeta )>n/2-N$, and is
also holomorphic in ${\mathbb C}-\frac{1}{2}(n-{\mathbb N})$.

Finally, one inverts $(I-E_\zeta )$ using analytic Fredholm theory.
Composition theorems of \cite{Mazzeo:1991a} show that if 
\[
(I+F_\zeta )=(I-E_\zeta )^{-1}, 
\]
then $\kappa (F_\zeta )$ also lies in $\eta ^\zeta (\eta ^{\prime }r)^\infty
C^\infty (X\times _0X)$. This in turn implies that $C_\zeta =(A_\zeta
+B_\zeta )F_\zeta $ belongs to $\beta ^{*}\left[ (xx^{\prime })^\zeta
C^\infty (\overline{X}\times \overline{X})\right] $, and therefore $\beta
^{*}G_\zeta $ has the claimed form.
\endproof

\begin{remark}
\label{rem.mapping}It follows from the form of the resolvent kernel that $
R_{V}(\zeta )$ is a continuous mapping from $\dot{C}^{\infty }(X)$ to $
C^{\infty }(X)$ when defined, and extends to a bounded
mapping from $x ^{N}L^{2}(X)$ to $x ^{-N}L^{2}(X)$ for $\Re (\zeta
)>n/2-N$.
\end{remark}

\begin{remark}
The Mazzeo-Melrose construction does not rule out the possibility of poles
at $\zeta \in \frac{1}{2}(n-{\mathbb N})$, possibly of infinite rank. If $(X,g)$
has constant negative curvature in a neighborhood of infinity, the operator $
A_{\zeta }+B_{\zeta }$ may be replaced by the model resolvent (the resolvent
of the Laplacian on the covering space ${\mathbb H}^{n+1}$), which is entire if 
$n$ is even and has finite-rank poles at $\zeta =-k\,$ if $n$ is odd (see
for example the explicit formulas in \cite{GZ:1995b}, section 2) . In either
case, these terms contain only poles of finite rank, and the last step of
the construction, involving the meromorphic Fredholm theorem, gives at most
poles with finite-rank residues. A detailed construction of the resolvent in
this case is given in \cite{GZ:1995b}, section 3. This observation plays a
crucial role in the proof of Theorem \ref{thm.hyperbolic}.
\end{remark}

Theorem \ref{thm.resolvent.ac} and standard arguments (see \cite{GZ:1997}
, Lemma 2.4) enable us to characterize the polar part of $R_{V}(\zeta )$ at
a resolvent resonance $\zeta _{0}\notin \frac{1}{2}(n-{\mathbb N})$. We will
view the meromorphically continued resolvent as a mapping from the space $
B_{0}=x ^{N}L^{2}(X)$ to $B_{1}=x ^{-N}L^{2}(X)$ as in Remark \ref
{rem.mapping}, where $N$ is chosen so that $\Re(\zeta_0)>n/2-N$. Introduce the nondegenerate form 
\[
\left\langle u,v\right\rangle =\int_{X}\,uv\,dx 
\]
(no complex conjugation) which can be used to pair elements in $B_{0}$ and $
B_{1}$. The resolvent operator is symmetric with respect to this form.

\begin{proposition}
\label{prop.resolvent.laurent}Let $\zeta _{0}\in {\mathbb C}-\frac{1}{2}(n-
{\mathbb N})$ be a pole of $R_{V}(\zeta )$. Then 
\[
R_{V}(\zeta )=\sum_{j=-k}^{-1}(\zeta (n-\zeta )-\zeta _{0}(n-\zeta
_{0}))^{j}\,A_{j}+H(\zeta ) 
\]
where $H(\zeta )$ is a holomorphic ${\mathcal L}(B_{0},B_{1})$-valued function
near $\zeta =\zeta _{0}$ and the $A_{j}$ are finite-rank operators in ${\mathcal 
L}(B_{0},B_{1})$ with 
\[
A_{-j}=(\Delta _{g}+V-\zeta _{0}(n-\zeta _{0}))^{j-1}A_{-1} 
\]
for $j\geq 2$. The operator $A_{-1}$ commutes with $\Delta _{g}+V-\zeta
_{0}(n-\zeta _{0})$. Moreover, there is a basis $\left\{ \psi _{i}\right\}
_{i=1}^{m_{\zeta _{0}}}$ for ${{\Ran}}\,(A_{-1})$ so that 
\[
A_{-1}f=\sum_{i}\left\langle f,\psi _{i}\right\rangle \,\psi _{i}, 
\]
and the operator $(\Delta _{g}+V-\zeta _{0}(n-\zeta _{0}))\,$ is represented
on ${{\Ran}}\,(A_{-1})$ by a matrix $M$ with $M^{k-1}\neq 0$ but $
M^{k}=0$.
\end{proposition}

\begin{remark}
If $(X,g)$ has constant curvature in a neighborhood of infinity, 
the same result holds for any resolvent resonance 
$\zeta _{0}$ including those $\zeta _{0}\in \frac{1}{2}(n-{\mathbb N})$.
\end{remark}

\section{The Scattering Operator for $P_V$}

\label{sec.scattering}

Let $S_{V}(\zeta )$ denote the scattering operator for $P_V=\Delta _{g}+V$. To
describe the scattering operator and its singularities, we recall its
connection with the resolvent kernel. First, we blow up $\partial \overline{X
}\times \partial \overline{X}$ along the diagonal $\Lambda _{\infty }$ to
obtain the space $\partial \overline{X}\times _{0}\partial \overline{X}$.
The map $\beta _{\partial }:\partial \overline{X}\times _{0}\partial 
\overline{X}\rightarrow \partial \overline{X}\times \partial \overline{X}$
is the `blow-down map' for this resolution. If $(y,y^{\prime })$ are local
coordinates for $\partial \overline{X}\times \partial \overline{X}$, $
r=\left| y-y^{\prime }\right| $, and $\theta =(y-y^{\prime })/r$, then $
(r,\theta ,y)$ give local coordinates for $\partial \overline{X}\times
_{0}\partial \overline{X}$, and $\beta _{\partial }(r,\theta
,y)=(y,y+r\theta )$. The kernel of the scattering operator is recovered as
an asymptotic limit of the resolvent kernel. Let $\kappa (A)$ denotes the
lift of the integral kernel of $A$ (with respect to the measure on $\partial 
\overline{X}\,$ induced by the metric $\left. h\right| _{\partial \overline{X
}}$) to $\partial \overline{X}\times _{0}\partial \overline{X}$. Then \cite
{Borthwick:1997,JS:1998} 
\begin{equation}
\kappa (S_{V}(\zeta ))=\left. \beta _{\partial }^{\ast }\left( (x x
^{\prime })^{-\zeta }G_{\zeta }\right) \right| _{T\cap B}  \label{eq.s.limit}
\end{equation}
where $T\cap B$ is the intersection of the top and bottom faces of $
\overline{X}\times _{0}\overline{X}$. From this formula and Theorem \ref
{thm.resolvent.ac}, we easily obtain:

\begin{theorem}
\label{thm.s.kernel}The decomposition 
\[
\kappa (S_{V}(\zeta ))=r^{-2\zeta }F_{\zeta }+\beta _{\partial }^{\ast
}(K_{\zeta }) 
\]
holds, where $F_{\zeta }$ and $K_{\zeta }$ are meromorphic maps respectively
into $C^{\infty }(\partial \overline{X}\times _{0}\partial \overline{X})$
and $C^{\infty }(\partial \overline{X}\times \partial \overline{X})$, and $
F_{\zeta }$ is holomorphic in ${\mathbb C}-\frac{1}{2}{\mathbb Z}$. At
poles $\zeta _{0}\notin {\mathbb C}-\frac{1}{2}{\mathbb Z}$, the kernel $
K_{\zeta }$ has polar part with finite Laurent series and coefficients in $
C^{\infty }(\partial \overline{X}\times \partial \overline{X})$. The set of
such poles is contained in the set of poles of $R_{V}(\zeta )$. \ For $\zeta
_{0}\in Z_{p}$ with $\zeta _{0}\notin \frac{1}{2}n+{\mathbb Z}$, $S_{V}(\zeta )$
has a semi-simple pole at $\zeta _{0}$ whose residue has rank $m_{\zeta
_{0}} $.
\end{theorem}

Note that the distribution $r^{-2\zeta }$ has poles at $\zeta \in \frac{1}{2}
n+{\mathbb N}$, giving rise to infinite-rank first-order poles of $S_{V}(\zeta )
$. The statement about poles at $\zeta _{0}\in Z_{p}$ follows from the facts
that the resolvent has a semi-simple pole at each such $\zeta _{0}$ and that
the  residue is a finite-rank projection onto eigenfunctions of the form $
x ^{\zeta _{0}}\psi $ for $\psi \in C^{\infty }(\overline{X})$ with $\left.
\psi \right| _{\partial \overline{X}}\neq 0$. Thus from (\ref{eq.s.limit}) it
follows that $S_{V}(\zeta )$ has a semi-simple pole with residue $
\sum_{i=1}^{m_{\zeta _{0}}}\left\langle \varphi _{i},\cdot \right\rangle
\varphi _{i}$ where $\varphi _{i}=\left. x ^{-\zeta _{0}}\psi _{i}\right|
_{\partial X}$.

It follows from Theorem \ref{thm.s.kernel} that the scattering operator may
be factored as 
\[
S_{V}(\zeta )=P(\zeta )(I+K_{V}(\zeta ))Q(\zeta )
\]
where $P(\zeta )$ and $Q(\zeta )$ are holomorphic families of invertible
elliptic operators in ${\mathbb C}-\frac{1}{2}{\mathbb Z}$, and $K_{V}(\zeta )$
is a meromorphic compact operator-valued function in ${\mathbb C}-\frac{1}{2}(n+
{\mathbb Z})$ with finite polar parts and finite-rank residues at each pole $
\zeta _{0}\in {\mathbb C}-\frac{1}{2}{\mathbb Z}$. For any 
$\zeta _{0}\notin \frac{1}{2}{\mathbb Z}$, 
\begin{eqnarray*}
\nu _{\zeta _{0}} &=&{{\Tr}}\left( \frac{1}{2\pi i}\int_{\gamma
_{\zeta _{0},\varepsilon }}S_{V}(\zeta )^{-1}S_{V}^{\prime }(\zeta )\,d\zeta
\right)  \\
&=&{{\Tr}}\left( \frac{1}{2\pi i}\int_{\gamma _{\zeta
_{0},\varepsilon }}(I+K_{V}(\zeta ))^{-1}K_{V}^{\prime }(\zeta )\,d\zeta
\right) ;
\end{eqnarray*}
since $\zeta \mapsto (I+K_{V}(\zeta ))$ is a meromorphic family of Fredholm
operators, it follows from \cite{GS:1972} that $\nu _{\zeta _{0}}$ is an
integer which may be thought of roughly as the number of zeros minus the
number of poles of $(I+K_{V}(\zeta ))$ at $\zeta =\zeta _{0}$. It may be
computed as follows. Near a pole 
$\zeta _{0}\notin \frac{1}{2}(n-{\mathbb N})$, the decomposition 
$$
(I+K_{V}(\zeta ))=
	E(\zeta )
	\left[\sum_{\ell }
		(\zeta -\zeta _{0})^{\nu _{\ell}}
		P_{\ell }+H(\zeta )\right]
	F(\zeta )
$$ 
for holomorphically invertible
operator-valued functions $E(\zeta )$ and $F(\zeta )$, finite-rank
projections $P_{\ell }$ with $P_\ell P_m=\delta_{\ell m}$, integers $\nu _{\ell }$, and a holomorphic
operator-valued function $H(\zeta )$; moreover a similar decomposition 
holds for $(I+K_{V}(\zeta ))^{-1}$ with the signs of $\nu _{\ell }$ 
reversed \cite{GS:1972}. Note that the sum over $\ell $ is finite and that
some $\nu _{\ell }$ may be negative. If 
$m_{\ell }={{{\Rank}(}}P_{\ell })$ then 
$\nu_{\zeta _{0}}=\sum_{\ell }\nu _{\ell }m_{\ell }$. From the functional
equation $S_{V}(\zeta )S_{V}(n-\zeta )=I$, it follows that there are
holomorphically invertible operator-valued functions $U(\zeta )$ and 
$V(\zeta )$ in $\zeta \notin \frac{1}{2}{\mathbb Z}$ with 
$$
(I+K_{V}(n-\zeta ))^{-1}=U\left( \zeta \right) (I+K_{V}(\zeta ))V(\zeta )
$$
so that the relationship 
\begin{equation}
\nu _{\zeta _{0}}=-\nu _{n-\zeta _{0}}
\end{equation}
holds.

Note that if $n-\zeta _{0}\in Z_{p}$, then the integer $\nu _{\zeta _{0}}$
will be determined both by the poles of $(I+K_{V}(\zeta ))$ and those of $
(I+K_{V}(n-\zeta ))$, which appear as zeros of $(I+K_{V}(\zeta ))$ at $\zeta
=\zeta _{0}$.

A careful analysis of the resolvent parametrix construction (see for example 
\cite{Borthwick:1997}) shows that the map 
\[
({\mathbb C}-\frac{1}{2}{\mathbb Z})\times \dot{C}^{\infty }(X)\ni (\zeta
,V)\longmapsto K_{V}(\zeta ) 
\]
is a continuous mapping away from poles of $K_{V}(\zeta )$.

Suppose now that $\zeta _{0}\in {\mathbb C}-\frac{1}{2}{\mathbb Z}$ is a 
{\em simple} resonance of $R_{V}(\zeta )$, $\zeta _{0}\notin Z_{p}$. The polar
part of $R_{V}(\zeta )$ is 
\[
(\zeta (n-\zeta )-\zeta _{0}(n-\zeta _{0}))^{-1}\left\langle \psi ,\cdot
\right\rangle \psi 
\]
where $\psi \in x ^{\zeta }C^{\infty }(\overline{X})$ solves the
eigenvalue equation\thinspace $(\Delta _{g}+V-\zeta _{0}(n-\zeta _{0}))\psi
=0$. We claim that $\varphi =\left. x ^{-\zeta }\psi \right| _{\partial 
\overline{X}}$ is a nonzero element of $C^{\infty }(\partial \overline{X})$.
If not, then introducing local coordinates $(x,y)$ on $X$ with $x=x $, we
have $\psi \in x^{\zeta +1}C^{\infty }$ and a power series argument using
the eigenvalue equation shows that, in fact, $\psi \in \dot{C}^{\infty }(X)$
, so that $\psi $ is an $L^{2}$-eigenfunction, a contradiction. It now
follows from (\ref{eq.s.limit}) that $S_{V}(\zeta )$ has a first-order pole
at $\zeta =\zeta _{0}$ with rank-one residue $(\varphi ,\cdot )\varphi $,
where $(\,\cdot \,,\,\cdot \,)$ is the real inner product on functions with
Riemannian measure induced by the metric $\left. h\right| _{\partial 
\overline{X}}$. From the holomorphic invertibility of $P(\zeta )$ and $
Q(\zeta )$ near $\zeta =\zeta _{0}$, it follows that $K(\zeta )$ has a
simple pole at $\zeta =\zeta _{0}$. Thus:

\begin{proposition}
\label{prop.simple}If $R_{V}(\zeta )$ has a simple resonance at 
$\zeta _{0}\in {\mathbb C}-\frac{1}{2}{\mathbb Z}$, then $
S_{V}(\zeta )$ also has a simple resonance at $\zeta _{0}$. If $R_V(\zeta)$ is
holomorphic at $\zeta_0 \in {\mathbb C}-\frac{1}{2}{\mathbb Z}$, then so
is $S_V(\zeta)$.
\end{proposition}

\begin{remark}
\label{rem.simple.cc}If $(X,g)$ has constant curvature near infinity,
Proposition \ref{prop.simple} holds for any simple resolvent resonance with 
$\zeta \not\in n/2-{\mathbb N}$ (i.e., including any simple poles at points 
$\zeta \in {\mathbb Z}$; later, we will analyze $\zeta \in n/2-{\mathbb N}$ separately).
The only detail to check is the
argument that shows that $\varphi =\left. x ^{-\zeta }\psi \right|
_{\partial \overline{X}}$ is a nonzero element of $C^{\infty }(\partial 
\overline{X})$. Since $\zeta _{0}$ and $(n-\zeta _{0})$ differ by an
integer, one can no longer show that $\psi \in \dot{C}^{\infty }(X)$ by a
formal power series argument, but one can show that if $\varphi =0$, then $
\psi \in x^{n}C^{\infty }(\overline{X})$, which is sufficient to show that $
\psi \in L^{2}(X)$ and derive a contradiction. 
\end{remark}

Note that, if {\em all} resonances of $R_{V}(\zeta )$ were simple,
Proposition \ref{prop.simple} and the conclusion of Theorem
\ref{thm.s.kernel}
that the set of scattering
poles is contained in the set of resolvent resonances would imply equality of
these two sets away from $\frac{1}{2}(n-{\mathbb N})$ and those $\zeta 
$ with $n-\zeta \in Z_{p}$. In the next two sections we shall show that 
all resonances are simple for `generic' $V$. The following continuity result 
for scattering poles will be useful.

\begin{proposition}
\label{prop.scatt.continuous}Let $\zeta _{0}\in {\mathcal S}$ with $\zeta
_{0}\in {\mathbb C}-\frac{1}{2}(n-{\mathbb N}{\bf )}$. There is a $\delta >0$ so
that for all $V\in \dot{C}^{\infty }(X)$ with $d(V,0)<\delta $ and that some $
\varepsilon >0$, the projection 
\[
\nu _{\zeta _{0}}(V)={\Tr}\left( \frac{1}{2\pi i}\int_{\gamma
_{\zeta _{0},\varepsilon }}S_{V}^{-1}(\zeta )S_{V}^{\prime }(\zeta )\,d\zeta
\right) 
\]
is continuous as a map from $\dot{C}^{\infty }(X)$ to ${\mathbb Z}$. In
particular, $\nu _{\zeta _{0}}(V)=\nu _{\zeta _{0}}(0)$ for such $V$.
\end{proposition}

\begin{proof}
This is a consequence of the continuity of the map $(\zeta ,V)\longmapsto
S_{V}(\zeta )$.
\end{proof}

\section{Perturbation Theory of Resonances}

\label{sec.perturb}

Now we apply Agmon's perturbation theory of resonances \cite{Agmon:1997} to
study the behavior of resonances under potential perturbations. For a fixed
asymptotically hyperbolic manifold $(X,g)$ we consider the family of
operators $P_{V}$ where $V$ ranges over complex-valued $\dot{C}
^{\infty }(X)$ functions. For any such $V$, Theorem \ref{thm.resolvent.ac}
guarantees that $R_{V}(\zeta )=(P_{V}-\zeta (n-\zeta ))^{-1}$ admits a
meromorphic continuation to any half-plane $\Re (\zeta )>n/2-N$, $N$ a
positive integer, as a mapping from $B_{0}=x ^{N}L^{2}(X)$ to $B_{1}=x
^{-N}L^{2}(X)$. We will set ${\mathcal R}_{V}(z)=(P_{V}-z)^{-1}$ with the
understanding that $z$ lies on the second sheet of the Riemann surface for
the inverse function of $f(\zeta )=\zeta (n-\zeta )$, so that ${\mathcal R}
_{V}(z)$ is the meromorphic continuation of the resolvent to the second
sheet. Using the fact that $P_{V}:\dot{C}^{\infty }(X)\rightarrow \dot{C}
^{\infty }(X)$ together with Theorem \ref{thm.resolvent.ac}, it is not
difficult to check that the operator $P_{V}$ satisfies the hypotheses of
Agmon's abstract theory.

To study the perturbation of a resonance $z_{0}$, Agmon introduces auxiliary
operators and Banach spaces associated to an open connected domain $\Delta $
containing $z_{0}$ with $C^{1}$ boundary $\Gamma $. Let $B_{\Gamma }$ be the
subset of $B_{1}$ consisting of functions of the form 
\[
u=f+\int_{\Gamma }{\mathcal R}_{V}(w)\,\Phi (w)\,dw 
\]
where $f\in B_{0}$ and $\Phi \in C(\Gamma ;B_{0})$, the continuous functions
on $\Gamma $ with values in $B_{0}$.  Finally, let 
$Y$ be the closed subspace of $B_0 \times C(\Gamma,B_0)$ consisting
of those $(g,\Phi )$ with 
\[
0=g+\int_{\Gamma }{\mathcal R}_{V}(w)\,\Phi (w)\,dw.
\]
The space $B_{\Gamma }$ is a Banach
space as the quotient of $B_{0}\times C(\Gamma ,B_{0})$ by the closed
subspace $Y$ when equipped with the quotient norm 
\[
\left\| u\right\| _{B_{\Gamma }}=\inf \left\{ \left\| \,f\,\right\|
_{B_{0}}+\left\| \Phi \right\| _{C(\Gamma ;B_{0})}:u=f+\int_{\Gamma }{\mathcal R}
_{V}(w)\,\Phi (w)\,dw\right\} . 
\]
The space $B_{\Gamma }$ satisfies $B_{0}\subset B_{\Gamma }\subset B_{1}$,
where the canonical injections are continuous.

The theory of \cite{Agmon:1997} implies that there is a closed operator $
P_{V}^{\Gamma }:{\mathcal D}(P_{V}^{\Gamma })\rightarrow B_{\Gamma }$ which is a
restriction of $P_{V}$ is a sense we will make precise, and whose
eigenvalues in $\Delta $ are exactly the resonances of ${\mathcal R}_{V}(z)$ in $
\Delta $. In fact, let $\overline{P}_{V}$ be the closure of $P_{V}$ as a
densely defined operator from $B_{1}$ to itself. Then $P_{V}^{\Gamma
}u=\overline{P}_{V}u$ for all $u\in {\mathcal D}(P_{V}^{\Gamma })$. The Laurent expansion
of ${\mathcal R}_{V}^{\Gamma }(z)=(P_{V}^{\Gamma }-z)^{-1}$ near a resonance $
z_{0}\in \Delta $ takes the form 
\[
\sum_{j=-k}^{-1}(z-z_{0})^{j}A_{j}^{\Gamma }+H^{\Gamma }(z) 
\]
where $H^{\Gamma }(z)$ is a holomorphic ${\mathcal L}(B_{\Gamma })$-valued
function in a neighborhood of $z_{0}$ and the $A_{j}^{\Gamma }$ are
finite-rank operators belonging to ${\mathcal L}(B_{\Gamma },{\mathcal D}
(P_{V}^{\Gamma }))$. For\thinspace $f\in B_{0}$ we have 
\[
A_{j}^{\Gamma }f=A_{j}f 
\]
where $A_{j}$ are the corresponding Laurent coefficients for ${\mathcal R}
_{V}(z) $.

Now fix $V\in \dot{C}^{\infty }(X)$ and set $P(t)=P_{V}+tW$ for another
potential $W\in \dot{C}^{\infty }(X)$. The operators $P(t)$ have resolvents
which admit meromorphic continuation to ${\mathcal L}(B_{0},B_{1})$-valued
meromorphic functions in $\Re (s)>n/2-N$ for any fixed $N>0$. Moreover, for $
t$ small and a fixed region $\Delta $, the spaces $B_{\Gamma }(t)$
corresponding to $P(t)$ are equal as sets and carry equivalent norms, and
the operators $P_{\Gamma }(t)$ form an analytic family of type (A)\ in the
sense of Kato \cite{Kato:1976}. Let ${\mathcal R}(t,z)=(P(t)-z)^{-1}$ and ${\mathcal 
R}^{\Gamma }(t,z)=(P^{\Gamma }(t)-z)^{-1}$. Theorem 7.7 of \cite{Agmon:1997}
shows that, for small $t$, the resolvents ${\mathcal R}^{\Gamma }(t,z)$ and $
{\mathcal R}(t,z)$ coincide on $B_{0}$, possess the same set of poles for each
fixed $t$, and ${{\Ran}(}A_{j}^{\Gamma }(t))={{\Ran}(}A_{j}(t))
$, where $A_{j}(t)$ and $A_{j}^{\Gamma }(t)$ are the respective Laurent
coefficients of ${\mathcal R}(t,z)$ and ${\mathcal R}^{\Gamma }(t,z)$ at a given
pole in $\Delta $.

\section{Generic Simplicity of Resonances}

\label{sec.generic}

For $L^{2}$ eigenvalues of the Laplacian and its perturbations, it has long
been known that `generic' potential perturbations split degenerate
eigenvalues so that a single eigenvalue of multiplicity $m$ becomes $m$
simple eigenvalues, localized near the original eigenvalue (see for example
\cite{Uhlenbeck:1976}, where generic simplicity is proved for the Laplacian
on compact manifolds, and Kato \cite{Kato:1976} for the background in
perturbation theory of linear operators; Uhlenbeck's methods adapt without
difficulty to eigenvalues below the continuous spectrum). The purpose of 
this section is to show that the same
is true of the resolvent resonances.

\begin{theorem}
\label{thm.simple}The set $E$ of potentials $V\in \dot{C}^{\infty }(X)$ for
which all eigenvalues of $\Delta _{g}+V$ and all resonances of $\Delta
_{g}+V $ in ${\mathbb C-}\frac{1}{2}(n-{\mathbb N})$ are simple is open and
dense in $\dot{C}^{\infty }(X)$.
\end{theorem}

We will follow rather closely the argument of \cite{KZ:1995} except that
Agmon's perturbation theory replaces the exterior complex scaling used
there. Since generic simplicity results for eigenvalues are well-known we
will only prove the genericity for the resonances.

For positive integers $N$ and real numbers $r>0$, we define 
\[
{\mathcal R}_{N}^{r}=\left\{ \zeta \in {\mathcal R}:\left| \zeta \right| <r,\,\,\,{
{\dist}}(\zeta ,\frac{1}{2}(n-{\mathbb N}))>1/N\right\} 
\]
and let 
\[
E_{N}^{r}=\left\{ V\in \dot{C}^{\infty }(X):
\mbox{each $\zeta$ in ${\mathcal R}_{N}^{r}$ is simple} \right\} .
\]
We set 
\[
E=\cap _{n=1}^{\infty }\cap _{N=1}^{\infty }E_{N}^{n},
\]
and we define 
\[
F=\dot{C}^{\infty }(X)\backslash E.
\]
We wish to show that $E$ is dense in $\dot{C}^{\infty }(X)$, i.e., that $F$
has empty interior. By the Baire category theorem, it suffices to show that $
F_{N}^{n}=\dot{C}^{\infty }(X)\backslash E_{N}^{n}$ is nowhere dense for
each $n$ and $N$. By the discreteness of the resonance set in ${\mathbb C}-
\frac{1}{2}(n-{\mathbb N})$, it suffices to show that for any $V\in F_{N}^{n}$,
any non-simple pole $\zeta _{0}$, and any $\varepsilon >0$, there is a $W$
with $x (W,0)<\varepsilon $ so that $V+W$ has only simple poles in a
neighborhood of $\zeta _{0}$. As in Section \ref{sec.perturb}, it will be convenient to work
with the meromorphically continued resolvent ${\mathcal R}(z)$, and we will
denote by $z_{0}$ the point on the second sheet of the Riemann surface for $
{\mathcal R}(z)$ corresponding to $\zeta _{0}$.

Consider the family of operators $P_{V+W}=\Delta _g+V+W$, where $d
(W,0)<\varepsilon _0$, and $\varepsilon _0>0$ is to be chosen (here 
$d(\cdot,\cdot)$ is the metric on $\dot C^\infty(X)$ defined in 
(\ref{eq.seminorm})). 
Let $\Gamma \,$ be a
contour enclosing $z_0$ and no other pole of ${\mathcal R}(z)$. 
For $\varepsilon _0$ small enough, the operators $P_{V+W}^\Gamma $ 
defined in Agmon's abstract
theory can be considered to act on a single Banach space $B_\Gamma $, and
the associated projection 
\[
\Pi _{V+W}^\Gamma =\frac 1{2\pi i}\int_\Gamma \,(P_{V+W}^\Gamma -w)^{-1}\,dw 
\]
is analytic in $W$ with $x (W,0)<\varepsilon _0$, and of constant rank,
say $m$. As in \cite{KZ:1995}, we note that either

\vskip 0.5cm

\begin{enumerate}
\item[(1)]  For each $\varepsilon >0$, there is a $W$ with $x
(W,0)<\varepsilon \,$ so that $P_{V+W}^\Gamma \,$ has at least two distinct
eigenvalues, or

\item[(2)]  There is an $\varepsilon >0$ so that for all $W$ with $x
(W,0)<\varepsilon $, $P_{V+W}^\Gamma $ has a single eigenvalue $z(W)$ and
there is an integer $k(W)$, $1\leq k(W)\leq m$, so that 
\[
(P_{V+W}^\Gamma -z(W))^{k(W)}\Pi _{V+W}^\Gamma =0\,\, \,\,(P_{V+W}^\Gamma
-z(W))^{k(W)-1}\Pi _{V+W}^\Gamma \neq 0\,. 
\]
\end{enumerate}
\vskip 0.5cm

If case (2) does not occur, we can split resonances repeatedly by small
perturbations. Thus we will suppose that case (2) does occur and obtain a
contradiction.

First, note that $k(W)$ is locally constant so by taking $\varepsilon _0$
small enough we may assume that $k(W)$ is constant for $W$ with $x
(W,0)<\varepsilon _0$. As in \cite{KZ:1995} we consider in turn the
possibilities $k(W)=1$ (the semi-simple case) and $k(W)\geq 2$.

First suppose that $k(W)=1$, that $z(W)\,$ is an eigenvalue of $
P_{V+W}^{\Gamma }\,$ of multiplicity $m$, and let $\left\{ \psi _{i}\right\}
_{i=1}^{m}$ be a basis for ${{\Ran}(}A_{-1}^{\Gamma })$, where $
A_{-1}^{\Gamma }$ occurs in the Laurent expansion for $(P_{V}^{\Gamma
}-z)^{-1}$ at $z=z_{0}$. The vectors $\left\{ \psi _{j}\right\} _{j=1}^{m}$
belong to $B_{\Gamma }\subset B_{1}$ and may be chosen to diagonalize $
A_{-1}^{\Gamma }$ as in Proposition \ref{prop.resolvent.laurent}. Let $
\left\{ f_{j}\right\} _{j=1}^{m}$ be a set of vectors in $B_{0}$ with $
\left\langle \psi _{i},f_{j}\right\rangle =\delta _{ij}$. Finally, for fixed 
$W$, let $L(t)=P_{V+tW}^{\Gamma }$, let $\Pi _{t}=\Pi _{V+tW}^{\Gamma }$,
let $\psi _{i}(t)=\Pi _{t}\psi _{i}$, and let $z(t)=z(tW)$. By
differentiating the eigenvalue equation 
\[
(L(t)-z(t))\psi _{i}(t)=0
\]
at $t=0$, we recover the identity 
\[
(W-z^{\prime }(0))\psi _{i}+(L(0)-z(0))\psi _{i}^{\prime }(0)=0.
\]
We now apply the projection $\Pi _{0}$ to both sides, pair with $f_{j}$, and
use the fact that $(L(0)-z(0))\Pi _{0}=\Pi _{0}(L(0)-z(0))=0$ to conclude
that 
\[
\left\langle f_{j},\Pi _{0}W\psi _{i}\right\rangle =z^{\prime }(0)\delta
_{ij}
\]
From the choice of $\left\{ f_{i}\right\} $ and the fact that $\Pi
_{0}=\sum_{i}\left\langle \psi _{i},\,\cdot \,\right\rangle \psi _{i}$ it
now follows that 
\[
\left\langle \psi _{i},W\psi _{j}\right\rangle =z^{\prime }(0)\delta _{ij}
\]
Since this must hold for any $W\in \dot{C}^{\infty }(X)$ (in particular for
all $W\in C_{0}^{\infty }(U)$ with $U$ an open subset of $X$), it follows
that at least one of the $\psi _{i}$ vanishes on $U$, and hence on $X$ by
unique continuation. This gives a contradiction.

Now suppose that $z(W)$ is not semi-simple, but that there is a fixed $k\geq
2$ so that 
\[
(L(t)-z(t))^{k}\Pi _{t}=0,\,\ \ \ (L(t)-z(t))^{k-1}\Pi _{t}\neq 0. 
\]
Choose a vector $h\in B$ with $\psi (t)=(L(t)-z(t))^{k-1}\Pi _{t}h\neq 0$,
so that $(L(t)-z(t))\psi (t)=0$. Let $\psi =\psi (0)$. A perturbation
calculation again leads to 
\begin{equation}
(W-z^{\prime }(0))\psi +(L-z(0))\psi ^{\prime }(0)=0.
\label{eq.nss.perturbed}
\end{equation}
If we apply the projection $\Pi _{0}$ to both sides of (\ref
{eq.nss.perturbed}) and pair with a vector $f\in B_{0}$ with $\Pi _{0}f=\psi 
$, \ we obtain 
\[
\left\langle \psi ,W\psi \right\rangle =z^{\prime }(0)\left\langle f,\Pi
_{0}\psi \right\rangle . 
\]
We have used the fact that $\Pi _{0}$ is symmetric with respect to the
pairing $\left\langle \,\cdot \,,\,\cdot \,\right\rangle $. To evaluate the
right-hand side, we use the fact that $L-z(0)$ preserves $B_{0}$ to write 
\begin{eqnarray*}
\left\langle f,\Pi _{0}\psi \right\rangle &=&\left\langle f,\Pi
_{0}(L-z(0))^{k-1}\Pi _{0}h\right\rangle \\
&=&\left\langle \psi ,\Pi _{0}(L-z(0)^{k-1})\Pi _{0}h\right\rangle \\
&=&\left\langle (L-z(0))^{k-1}\Pi _{0}h,(L-z(0)^{k-1})\Pi _{0}h\right\rangle
\\
&=&\left\langle (L-z(0))^{k-2}\Pi _{0}h,(L-z(0)^{k})\Pi _{0}h\right\rangle \\
&=&0
\end{eqnarray*}
so that $\left\langle \psi ,W\psi \right\rangle =0$ for all $W\in \dot{C}
^{\infty }(X)$. It follows that $\psi $ vanishes on $X$, a contradiction.

We have now shown that for each $\varepsilon >0$, there is a $W$ with $d
(W,0)<\varepsilon \,$ so that $P_{V+W}^{\Gamma }\,$ has at least two
distinct eigenvalues. It follows that any resonance can be split by an
arbitrarily small perturbation $W\in \dot{C}^{\infty }(X)$. This fact implies
that the set $E$ of potentials $V$ for which $\Delta _{g}+V$ has only simple
resonances in ${\mathbb C-}\frac{1}{2}(n-{\mathbb N)}$ is open and dense in $\dot{C
}^{\infty }(X)$, and Theorem \ref{thm.simple} is proved.

In case $(X,g)$ has constant curvature near infinity, this result can be
improved if we work with the class $C_{0}^{\infty }(U)$ for a fixed open
subset $U$ of $X$. In this case the methods of \cite{GZ:1995b} can be used
to show that the resolvent of $\Delta _{g}+V$ \ has a meromorphic
continuation with only finite-rank poles, including any poles at $\zeta
_{0}\in \frac{1}{2}(n-{\mathbb N})$. One can then apply the above arguments
without essential changes to prove:

\begin{theorem}
\label{thm.simple.hyperbolic}Let $U$ be a fixed open subset of $X$ with
compact closure. The set $
E$ of potentials $V\in C_{0}^{\infty }(U)$ for which all eigenvalues and all
resonances of $\Delta _{g}+V$ are simple is open and dense in $C_{0}^{\infty
}(U)$.
\end{theorem}

\section{Resolvent Resonances and Scattering Poles}

\label{sec.equality}

Finally, we give the proofs of Theorems \ref{thm.essential} and \ref
{thm.hyperbolic}.

To prove Theorem \ref{thm.essential}, we choose a $\dot{C}^{\infty }(X)$
potential $V$ so that all of the eigenvalues and resonances of $\Delta _{g}+V
$ are simple. We further choose $V$ small enough that, for a given point $(\zeta
_{0},m_{\zeta _{0}})\in {\mathcal R}$ with $\zeta _{0}\notin {\mathbb C-}\frac{1}{2}
(n-{\mathbb N)}$, some $\varepsilon >0$, and any $t\in (0,2)$, no resonances of 
$\Delta _{g}+tV$ cross the circle $\gamma _{\zeta _{0},\varepsilon }$ of
radius $\varepsilon $ about $\zeta _{0}$, and the projection 
\[
\Pi _{tV}=\frac{1}{2\pi i}\int_{\gamma _{\zeta _{0},\varepsilon }}(2\zeta
-n)(\Delta _{g}+tV-\zeta (n-\zeta ))^{-1}\,d\zeta 
\]
(as well as its analogue for $n-\zeta _{0}$ if $\zeta _{0}\in Z_{p}$) is
continuous in $t$. It follows from Kato-Rellich perturbation theory for
small $t$, the rank of $\Pi _{tV}$ is continuous, so that 
\[
m(t)={{\rank}(}\Pi _{tV})=m_{\zeta _{0}}
\]
is constant for $t$ small. On the other hand, for $t\neq 0$ and small, the
resonances of $\Delta _{g}+V$ are simple, and the scattering operator has $
m_{\zeta _{0}}$ simple poles near $\zeta _{0}$ by Proposition \ref
{prop.simple}, and $m_{n-\zeta _{0}}$ simple zeros by the remarks following 
(\ref{eq.s.mult}) and the fact that an eigenvalue of multiplicity $m_{n-\zeta
_{0}}$ splits into $m_{n-\zeta _{0}}$ simple eigenvalues near $n-\zeta _{0}$
under perturbation. Finally, by Proposition \ref{prop.scatt.continuous}, $
m_{\zeta _{0}}-m_{n-\zeta _{0}}=\nu _{\zeta _{0}}$, where $\nu _{\zeta _{0}}$
is the multiplicity of the scattering pole $\zeta _{0}$ of $S(\zeta )$, the
scattering operator for $\Delta _{g}$. Hence $\zeta _{0}\in {\mathcal S}$, and $
m_{\zeta _{0}}-m_{n-\zeta _{0}}=\nu _{\zeta _{0}}$. Theorem \ref
{thm.essential} is proved.

For poles $\zeta _{0}\in {\mathbb C-}\frac{1}{2}(n-{\mathbb N)}$, the proof breaks
down for several reasons: 
(i)  the resolvent may have infinite rank poles at these points, 
(ii) the scattering operator may have infinite rank zeros at the points
$\zeta \in n/2-{\mathbb N}$, and 
(iii) the proof of
Proposition \ref{prop.simple} breaks down. If $(X,g)$ has constant curvature
in a neighborhood of infinity, the dimension of $X$ is even, and the potential 
perturbation $V$ is compactly supported in $X$,  we can show that (i) any resolvent
resonances for $P_V$ occuring at points $\zeta \not\in n/2-{\mathbb N}$ are 
semi-simple, and so can be perturbed to simple resonances under small 
potential perturbations, and (ii)
the resolvent is actually holomorphic at $\zeta\in n/2-{\mathbb N}$, 
and the scattering operator $S_V(\zeta)$ vanishes at these points. 
These observations will enable us to
prove Theorem 1.2.

Following \cite{PP:1998}, we define a renormalized scattering operator by
setting 
\begin{equation}
S_{r,V}(\zeta )=2^{2\zeta -n}\frac{\Gamma (\zeta -n/2)}{\Gamma
(n/2-\zeta )}S_{V}(\zeta ).  \label{eq.new.scattering}
\end{equation}
It follows from \cite{JS:1998} that $S_{r,V}(\zeta )$ admits the
factorization 
\[
S_{r,V}(\zeta )=P(\zeta )(I+K_{V}(\zeta ))P(\zeta )
\]
where now $P(\zeta )$ is a holomorphic family of invertible operators with $
P(\zeta )P(n-\zeta )=I$, and $K_{V}(\zeta )$ is a meromorphic family of
compact operators with finite polar parts whose Laurent coefficients are
finite-rank operators. Moreover the map $(\zeta ,V)\longmapsto K_{V}(\zeta )$
is continuous away from poles of $K_{V}(\zeta )$. We now define ${\mathcal S}$
to be the set of poles of $(I+K_{0}(\zeta ))$ with multiplicities 
\[
\nu _{\zeta _{0}}=
{\Tr}\left( \frac{1}{2\pi i}
\int_{\gamma _{\zeta _{0}},\varepsilon }
	S_{r,0}(\zeta )^{-1} S_{r,0}^{\prime }(\zeta)\,d\zeta \right) .
\]
This definition coincides with the previous definition for all $\zeta
_{0}\in {\mathbb C}-\frac{1}{2}(n-{\mathbb N})$ and is well-defined for any 
$\zeta$. Note that the distribution kernel of $
S_{r,V}(\zeta )$ is recovered from the Schwarz kernel of the resolvent
by (compare (\ref{eq.s.limit})) 
\[
\kappa (S_{V}(\zeta ))=2^{2\zeta -n}\frac{\Gamma (\zeta -n/2)}{\Gamma
(n/2-\zeta )}\left. \beta ^{\ast }\left( (x x ^{\prime })^{-\zeta
}G_{\zeta }\right) \right| _{T\cap B}
\]
so that for $\zeta _{0}\not\in \frac{1}{2}n-{\mathbb N}$, Remark \ref{rem.simple.cc}
shows that if $R_{V}(\zeta )$ has a simple
pole at $\zeta _{0}$, then $S_{V}(\zeta )$ does also.
It remains to analyze what happens for
the points $\zeta_0 \in \frac{1}{2}n-{\mathbb N}$ and it is
here that our assumption on even dimension enters. In this case, we can
prove:

\begin{lemma} Suppose that $X$ has constant curvature in a neighborhood of
infinity, $V$ is compactly supported, and $\dim X$ is even. Then the resolvent
$R_V(\zeta)$ is holomorphic at $\zeta \in \frac{1}{2}n-{\mathbb N}$, and $S_V(\zeta)$
vanishes at these points.
\end{lemma}

\begin{proof}
This proof is inspired by the proof of \cite{GZ:1997}, Lemma 2.8, and
generalizes its argument.  Let $G_{0,\zeta}(w,w')$ be the integral kernel 
of $R_0(\zeta)$ with respect to Riemannian measure on ${\mathbb H}^{n+1}$, 
where $w=(y,x)$, $w'=(y',x')$ for $y,y;' \in {\mathbb R}^n$ and $x,x'>0$
(upper half-space model).  Then 
$G_{0,\frac{1}{2}n - k} \in 
(xx')^{\frac{1}{2}n+k}C^\infty({\mathbb H}^{n+1}\times{\mathbb
H}^{n+1} \setminus \Delta)$ where $\Delta$ denotes the diagonal. This may be
checked by explicit computation of the resolvent or by using the formula
$$ G_{0,\zeta}(w,w')-G_{0,n-\zeta}(w,w')=
\int_{{\mathbb R}^n} e_{0,\zeta}(w,y'')e_{0,n-\zeta}(w,y'') \, dy''
$$
where 
$$e_{0,\zeta}(w,y'')=
	\pi^{-n/2}\frac{\Gamma(\zeta)}{\Gamma(\zeta-n/2)}
	\frac{x^\zeta}{(|y-y''|^2+x^2)^\zeta}.
$$
The point is that $e_{0,\zeta}$ vanishes at $s \in \frac{1}{2}n-{\mathbb
N}$ if $\dim X=n+1$ is even.  It follows
from the iterative parametrix construction in \cite{GZ:1995b} that the 
meromorphically continued resolvent kernel of $R(\zeta)$ has no pole at
$\zeta=n/2-k$, since its residue would correspond to a solution of the
eigenvalue equation with $y^{n/2+k}$ (hence square-integrable) behavior at
infinity. It further follows that the resolvent kernel belongs to
$(xx')^{n/2+k}C^\infty(X \times X )$ away from the diagonal,
so that the Schwartz kernel of the scattering operator at $\zeta=n/2-k$,
obtained as a boundary value of the resolvent rescaled by
$(xx')^{\frac{1}{2}(k-n)}$, is zero, and hence the scattering operator is
zero.
\end{proof}

From the continuity of the map $(\zeta ,V)\longmapsto K_{V}(\zeta )$ it
follows that the integer-valued function 
\[
\nu _{\zeta _{0}}(V)={\Tr}\left( \frac{1}{2\pi i}\int_{\gamma
_{\zeta _{0}},\varepsilon }S_{r,V}(\zeta )^{-1}S_{r,V}^{\prime
}(\zeta )\,d\zeta \right) 
\]
is constant for some $\delta >0$ and all $V\in \dot{C}^{\infty }(X)$ with $
d(0,V)<\delta $. The continuity of the projection $\Pi _{V}$ for resolvent
poles can also be established for poles in $\frac{1}{2}(n-{\mathbb N)}$ since
any such poles are known to have finite rank. We are now able to argue as
before to complete the proof.

\bibliographystyle{amsplain}

\end{document}